\documentclass[12pt]{amsart}

\begin{document}
\title[Moore-Penrose Inverse ]{ On the Moore-Penrose
inverse in $C^*$-algebras }
\author[Enrico Boasso]{Enrico Boasso}
\begin{abstract}In this article, two results regarding the Moore-Penrose inverse in the 
frame of $C^*$-algebras are considered. In first place, a characterization of the so-called
reverse order law is given, which provides a solution of a problem posed by
M. Mbekhta. On the other hand, Moore-Penrose 
hermitian elements, that is $C^*$-algebra elements which coincide
with their Moore-Penrose inverse, are intro-
duced and studied. In fact,
these elements will be fully characterized both in the Hilbert space
and in the  $C^*$-algebra setting. Furthermore, it will be proved that
an element is normal and Moore-Penrose hermitian if and only if
it is a hermitian partial isometry. 

\end{abstract}
\maketitle \noindent\bf{1. Introduction}\rm \vskip.5cm

\indent Given an unitary ring $A$, an element $a\in A$ will be
called \it regular\rm, if it has a \it generalized inverse\rm, also called
\it pseudo-inverse, \rm in A,
that is if there exists $a{'}\in A$ for which
$$
a=aa{'}a.
$$

\noindent It is clear that in this case $aa{'}$ and $a{'}a$ are idempotents of
$A$.\par

\indent In addition, a generalized inverse $a{'}$ of a regular
element $a\in A$ will be called \it normalized\rm, if $a{'}$ is
regular and $a$ is a pseudo-inverse of $a{'}$, that is if
$$
a=aa{'}a,\hskip1cm a{'}=a{'}aa{'}.
$$
\indent In the presence of an involution $*\colon A\to A$, it is
also possible to enquire if the idempotents  $aa{'}$ and $a{'}a$ are
self-adjoint, equivalently whether or not
$$(aa{'})^*=aa{'},\hskip1cm (a{'}a)^*=a{'}a.$$
In this case $a{'}$ is called the \it Moore-Penrose inverse \rm of
$a$, and it is denoted by  $a^{\dag}$, see [16], where this
concept was introduced for matrices, and the related works [10], [11],
and [14].\par \indent In [10] it was proved that each regular element
$a$ in a $C^*$-algebra $A$ has a Moore-Penrose inverse, which in
addition is unique. Consequently, the Moore-Penrose inverse of a
regular element $a\in A$ is the unique solution $x\in A$ to the
following set of equations:
$$
a=axa,\hskip.5cm x=xax,\hskip.5cm (ax)^*=ax,\hskip.5cm (xa)^*=xa.
$$
\indent According to the uniqueness of the Moore-Penrose inverse of
a regu- lar element $a$, $a^*$ also has a Moore-Penrose inverse and
$$
(a^*)^{\dag}=(a^{\dag})^*.
$$
\indent Furthermore, according to the above equations, if $a$ is a
regular ele-ment, then $a^{\dag}$ also is and
$$
(a^{\dag})^{\dag}=a.
$$
\indent The so-called reverse order law is one of the most important 
properties of the Moore-Penrose inverse
that have been deeply studied,
that is under what condition the equation 
$$
(ab)^{\dag}=b^{\dag}a^{\dag}
$$
holds.\par
\indent In the well-known article [7], T. N. E. Greville proved that the following facts
are equivalent:\par
\vskip.15cm i-\hskip.2cm $(ab)^{\dag}=b^{\dag}a^{\dag}$,\par
ii- \hskip.1cm $a^{\dag}abb^*a^*=bb^*a^*$ and $bb^{\dag}a^*ab = a^*ab$,\par
iii- $a^{\dag}a$ commutes with $bb^*$ and $a^*a$ with $bb^{\dag}$,\par
iv- $a^{\dag}abb^*a^*abb^{\dag} = bb^*a^*a$,\par
v- \hskip.1cm $a^{\dag}ab=b(ab)^{\dag}ab$ and $bb^{\dag}a^*= a^*ab(ab)^{\dag}$,\par
 \vskip.15cm
\noindent where $a$ and $b$ are two matrices. However, it is worth
noticing that the proofs in [7]
are also valid in the more general context of $C^*$-algebras. \par

\indent The key results of [7] were extended in some works devoted to generalized 
inverses of matrices, see for example [2], [3], and [18]. As
regard Hilbert space operators, in [4] R. Bouldin gave a
characterization in terms of invariant subspaces, which was
refined in [12] and [17]. Observe that the main result in [4],
Theorem 3.1, is equivalent to the generalization of Theorem 2 of [7], 
the above mentioned condition iii, to Hilbert space operators, see
Remark 3.2 of [4].\par   
 
\indent On the other hand, in the work [14] M. Mbekhta studied the
reverse order law for generalized inverses in the
frame of $C^*$-algebras. In fact, given two
regular elements $a$ and $b$ in a $C^*$-algebra $A$,
it was proved that the following statements are equivalent:\par
\vskip.15cm \it i\rm-  \hskip.2cm $b{'}a{'}$ is a generalized inverse of
$ab$,\par \it ii\rm- \hskip.1cm $a(pq-qp)b=0$,\par \it iii\rm- $qp$  is an
idempotent,\par\vskip.15cm 
\noindent where $a{'}$ and $b{'}$ are generalized
inverses of $a$ and $b$ respectively, $p=bb{'}$ and $q=a{'}a$, see
Theorem 3.1 of [14]. 
Naturally, this characterization remains true in Banach algebras,
in fact in a ring. Furthermore, in [5] R. Bouldin proved the
same characterization for Banach space operators.\par
\indent In addition, in [14] M. Mbekhta posed the problem of finding
necessary and sufficient conditions, analogues to the ones of 
Theorem 3.1 of [14], which ensures that
$$ (ab)^{\dag}= b^{\dag}a^{\dag},
$$
\noindent for $a$ and $b$ in a $C^*$-algebra $A$.\par

\indent In the work [13] it was claimed that the question
of M. Mbekhta in [14] was solved. However, 
the answer to this problem, Theorem 5 of [13],
not only does not provide conditions analogues
to the one of Theorem 3.1 in [14], but also
it consists in the formulation of the well-known Theorems 1 and 2
of [7] in $C^*$-algebras, the above reviewed conditions i, ii, and iii,
whose proofs are also valid in $C^*$-algebras.\par

\indent The first and main objective of the
present work consists in solving the problem posed by 
M. Mbekhta, that is to give a
characterization of the reverse order law for the Moore-Penrose
inverse in $C^*$-algebras which is analogue to the one of Theorem 3.1 of [14]. 
Due to the fact that the Moore-Penrose inverse is determined by
four equations instead of one, and that it
involves not only the product but also the involution,
several modifications must be made, however
the form of M. Mbekhta's characterization is preserved. What is
more, in section 3 four equivalent characterizations with this
characteristic will be proved. To this end it will be
necessary to reformulate the equations that define the Moore-Penrose
inverse of a regular element, which will be done in section 2
following an argument in [16].\par 

\indent On the other hand, given a regular element $a$ in a $C^*$-algebra
$A$, according to a general argument, or even as an application of the
results of section 3, it is easy to prove that $(aa^{\dag})^{\dag}=
aa^{\dag}$ and $(a^{\dag}a)^{\dag}=a^{\dag}a$. Now well, since the
Moore-Penrose inverse is a particular generalized inverse, it can be
thought of a sort of inverse, however, these two identities also
suggest that the Moore-Penrose inverse has properties that are
similar to the ones of the involution of the algebra. This
observation has led to the second objective of this work, namely,
the study of the regular elements $a\in A$ for which $a^{\dag}=a$.
These elements will be called Moore-Penrose hermitian, and its
basic properties will be studied in section 4. Furthermore, in
section 5 Moore-Penrose hermitian elements will be fully
characterize both in the Hilbert space and in the $C^*$-algebra
setting. In addition, it will be also proved that $a\in A$ is a
normal Moore-Penrose hermitian element if and only if it is a
hermitian partial isometry.\par \indent This article was written
during a research visit to the Abdus
Salam International Centre for Theoretical Physics, thanks to the
Research Fellowships 2005 Programm. The author wishes to express his
indebtedness to the authorities of the Mathematics Section of the
ICTP. In fact, the stimulating atmosphere and the warm hospitality
of the mentioned centre were two extraordinary helps to the research
work of the author. \par \indent This work was also supported by
UBACyT and CONICET.\par \vskip.5cm
 \bf{2. Equivalent formulations of the Moore-Penrose inverse}\rm \vskip.5cm
\indent Consider a $C^*$-algebra $A$, and $a\in A$ a regular 
element. In this section
several equivalent formulations of the equations defining the
Moore-Penrose inverse of $a$ will be considered. These
formulations will be central in the proof of the characterizations
of the next section. In addition, the argument in Proposition 2.1
will follow ideas of Theorem 1 of [16].\par
\newtheorem*{prop2.1}{Proposition 2.1}
\begin{prop2.1} Consider a  $C^*$-algebra $A$,
and two elements in $A$, $a$ and $x$. Then,\par \indent i-
\hskip.15cm the equations $a=axa$ and $(ax)^*=ax$ are equivalent to
$ a=$\par \hskip.8cm $x^*a^*a$.\par \indent ii- \hskip.1cm the
equations $a=axa$ and $(xa)^*=xa$ are equivalent to $ a=$
\par\hskip.8cm $aa^*x^*$.\par
 \indent iii-
 the
equations $x=xax$ and $(ax)^*=ax$ are equivalent to $ x=$
\par\hskip.1cm $xx^*a^*$.\par \indent \hskip.1cm iv- the equations
$x=xax$ and $(xa)^*=xa$ are equivalent to $ x=$ \par\hskip.8cm
$a^*x^*x$.\par
\end{prop2.1}
\begin{proof}
\indent The third equivalence was proved in Theorem 1 of [16]. The
other three statements can be proved in a similar way.\par
\end{proof}
\indent As a consequence, the following equivalent conditions are
obtained.\par
\newtheorem*{prop2.2}{Proposition 2.2}
\begin{prop2.2} Consider a $C^*$-algebra $A$ and $a\in A$. Then the
following statements are equivalent: \par \indent i- \hskip.2cm
$x\in A$ is the Moore-Penrose inverse of $a$,\par \indent ii-
\hskip.1cm $ a=x^*a^*a$ and $ x=a^*x^*x$,\par \indent iii- $
a=aa^*x^*$ and $x=xx^*a^*$.
\end{prop2.2}
\begin{proof}
\indent Is is a consequence of Proposition 2.1 and the equations
defining the Moore-Penrose inverse.\par
\end{proof}
\newtheorem*{rem2.3}{Remark 2.3}
\begin{rem2.3}\rm Consider a $C^*$-algebra $A$, $a\in A$ a regular
element of $A$, and $x=a^{\dag}$. Then, according to Proposition 2.2
and to the fact that $a^*$ is also regular and $(a^*)^{\dag}$ =
$(a^{\dag})^*$, the following statements are equivalent:\par \indent
\it{i-} \hskip.2cm $x\in A$ is the Moore-Penrose inverse of $a$,\par
\indent  \it{ii-} \hskip.15cm $a^*=a^*ax$ and $x^*=x^*xa$,\par
\indent {iii-} $ a^*=xaa^*$ and $ x^*=axx^*$.
\end{rem2.3}
\indent Next follows the equivalent formulations of the
Moore-Penrose inverse that will be central in the next section.\par
\newtheorem*{prop2.4}{Proposition 2.4}
\begin{prop2.4}Consider a $C^*$-algebra $A$ and $a\in A$. Then the following
statements are equivalent:\par \indent i- \hskip.2cm $x\in A$ is the
Moore-Penrose inverse of $a$,\par \indent ii- \hskip.1cm$a^*=xaa^*$
and $x=xx^*a^*$,\par \indent iii- $a=aa^*x^*$ and
$x^*=axx^*$,\par\indent iv- \hskip.1cm $a^*=a^*ax$ and
$x=a^*x^*x$,\par \indent v- \hskip.2cm $a=x^*a^*a$ and
$x^*=x^*xa$.\par
\end{prop2.4}
\begin{proof} It is a consequence of Proposition 2.2 and Remark 2.3.
\end{proof}
\vskip.5cm \noindent \bf{3. The reverse order law for the
Moore-Penrose inverse}\rm \vskip.5cm

\indent In this section the relationship between the product and the
Moore-Penrose inverse will be studied. In fact, four equivalent
characterizations of the so-called reverse order law for the
Moore-Penrose inverse will be proved. These characterizations are
analogue to the one given in Theorem 3.1 of [14] for the generalized
inverse of the product of two $C^*$-algebra elements. The results of
this section provide an answer to a question posed by M. Mbekhta in
[14].
\par

\newtheorem*{theo3.1}{Theorem 3.1}
\begin{theo3.1} Consider a $C^*$-algebra $A$, and two regular elements of $A$, $a$ and $b$,
such that $ab$ is also regular. Define $p$ = $bb^{\dag}$, $q$ =
$a^{\dag}a^{\dag^*}$, $r$ = $bb^*$ and $s$ = $a^{\dag}a$. Then,
the following statements are equivalent:\par \indent i- \hskip.3cm
$(ab)^{\dag} =b^{\dag}a^{\dag}$,\par \indent ii- \hskip.2cm
$a(pq-qp)b^{\dag^*}=0$, and $a(rs-sr)b^{\dag^*}=0$,\par \indent iii-
$spqp=qp$, and $srsp=sr$.\par
\end{theo3.1}
\begin{proof}
\indent First of all, observe that $p$, $q$, $r$ and $s$ are hermitian elements of
$A$.\par
\indent Consider $a^{\dag}$, $b^{\dag}$ and $(ab)^{\dag}$, the
Moore-Penrose inverses of $a$, $b$ and $ab$ respectively. According
to the third statement of Proposition 2.4, the following equations
hold:
\begin{align*}
&a=aa^*a^{\dag^*},  \hskip.5cmb=bb^*b^{\dag^*},\hskip.6cm ab= ab(ab)^*(ab)^{\dag^*},\\
&a^{\dag^*}=aa^{\dag}a^{\dag^*}, \hskip.3cm
b^{\dag^*}=bb^{\dag}b^{\dag^*},\hskip.2cm
(ab)^{\dag^*}=ab(ab)^{\dag}(ab)^{\dag^*}.
\end{align*}
\indent Furthermore, note that according again to the third
statement of Proposition 2.4, $$ a=as,\hskip.4cm
a^{\dag^*}=aq,\hskip.4cm b=rb^{\dag^*}, \hskip.4cm
b^{\dag^*}=pb^{\dag^*}. $$ \indent Now suppose that $(ab)^{\dag}$ =
$b^{\dag}a^{\dag}$. Then, since $(ab)^*$ = $b^*a^*$ and
$(ab)^{\dag^*}$ = $(b^{\dag}a^{\dag})^*$ = $a^{\dag^*}b^{\dag^*}$,
it is clear that
$$
ab= abb^*a^*a^{\dag^*}b^{\dag^*},\hskip1cm
a^{\dag^*}b^{\dag^*}=abb^{\dag}a^{\dag}a^{\dag^*}b^{\dag^*},
$$
which is equivalent to
$$
asrb^{\dag^*}=arsb^{\dag^*},\hskip1.5cm aqpb^{\dag^*}=
apqb^{\dag^*},
$$
which in turn is equivalent to the following identities:
$$
a(pq-qp)b^{\dag^*}=0, \hskip2cm a(rs-sr)b^{\dag^*}=0.
$$
\indent Next suppose that the second statement of the theorem holds.
Then, it is clear that
$$
a^{\dag}apqb^{\dag^*}b^*= a^{\dag}aqpb^{\dag^*}b^*,\hskip1cm
a^{\dag}arsb^{\dag^*}b^*= a^{\dag}asrb^{\dag^*}b^*.
$$
\indent However, according again to the third statement of
Proposition 2.4, and to the fact that $s=s^*$ and $p=p^*$, these equations can be rewritten as
$$
spqp=a^{\dag}(aa^{\dag}a^{\dag^*})b(b^{\dag}b^{\dag^*}b^*)=a^{\dag}a^{\dag^*}bb^{\dag}=qp,
$$
$$
srsp=(a^{\dag}aa^*)a^{\dag^*}(bb^*b^{\dag^*})b^*=a^*a^{\dag^*}bb^*=sr.
$$
\indent Finally suppose that the third statement of the theorem
holds. Then, since $p=p^*$ and $s=s^*$, it is clear that
$$
a^{\dag}abb^{\dag}a^{\dag}a^{\dag^*}b^{\dag^*}b^*=a^{\dag}a^{\dag^*}bb^{\dag},
$$
$$
a^{\dag}abb^*a^*a^{\dag^*}b^{\dag^*}b^*=a^*a^{\dag^*}bb^*,
$$
which implies that
$$
(aa^{\dag}a)bb^{\dag}a^{\dag}a^{\dag^*}(b^{\dag^*}b^*b^{\dag^*})=(aa^{\dag}a^{\dag^*})(bb^{\dag}b^{\dag^*}),
$$
$$
(aa^{\dag}a)bb^*a^*a^{\dag^*}(b^{\dag^*}b^*b^{\dag^*})=(aa^*a^{\dag^*})(bb^*b^{\dag^*}).
$$
However, according to the third statement of Proposition 2.4 and to
the fact that $a^{\dag}$ and $b^{\dag}$ are the Moore-Penrose
inverse of $a$ and $b$ respectively, the previous equations are
equivalent to
$$
ab(b^{\dag}a^{\dag})(b^{\dag}a^{\dag})^* = (b^{\dag}a^{\dag})^*,
$$
$$
ab(ab)^*(b^{\dag}a^{\dag})^*=ab,
$$
which, according again to the third statement of Proposition 2.4
implies that
$$
(ab)^{\dag}=b^{\dag}a^{\dag}.
$$
\end{proof}

\indent Note that in a $C^*$-algebra, under the same conditions
of Theorem 3.1, when instead of generalized inverses Moore-Penrose
inverses are considered, the characterization of Theorem 3.1 in [14]
determines if $b^{\dag}a^{\dag}$ is a normalized generalized inverse
of $ab$. However, in order to characterize the reverse order law for
the Moore-Penrose inverse, another equation is necessary
as well as new elements must be introduced. \par

\indent The next three theorems provide characterizations
which are equivalent to the one in Theorem 3.1. However,
for sake of completeness they are included.\par 

\newtheorem*{theo3.2}{Theorem 3.2}
\begin{theo3.2} Under the same conditions and notations of Theorem 3.1,
the following statements are equivalent:\par \indent i- \hskip.3cm
$(ab)^{\dag} =b^{\dag}a^{\dag}$,\par \indent ii- \hskip.2cm
$b^{\dag}(qp-pq)a^*=0$, and $b^{\dag}(sr-rs)a^*=0$,\par \indent
iii- $pqps=pq$, and $psrs=rs$.\par
\end{theo3.2}
\begin{proof}
\indent The proof is similar to the one of Theorem 3.1. However,
instead of the third statement of Proposition 2.4, the second
statement of the mentioned proposition must be used.\par
\end{proof}

\newtheorem*{theo3.3}{Theorem 3.3}
\begin{theo3.3} Under the same conditions and notations of Theorem 3.1, 
the following statements are equivalent:\par \indent i- \hskip.3cm
$(ab)^{\dag} =b^{\dag}a^{\dag}$,\par \indent ii- \hskip.2cm
$b^*(q^{\dag}p-pq^{\dag})a^{\dag}=0$, and
$b^*(sr^{\dag}-r^{\dag}s)a^{\dag}=0$,\par \indent iii-
$pq^{\dag}ps=pq^{\dag}$, and
$psr^{\dag}s=r^{\dag}s$.\par
\end{theo3.3}
\begin{proof}
\indent The proof is similar to the one of Theorem 3.1. However,
instead of the third statement of Proposition 2.4, the forth
statement of the mentioned proposition must be used. In addition, in
order to compute $q^{\dag}$ and $r^{\dag}$, Theorem 7 of [10] must be
considered.\par
\end{proof}

\newtheorem*{theo3.4}{Theorem 3.4}
\begin{theo3.4} Under the same conditions and notations of Theorem 3.1,
the following statements are equivalent:\par
\indent i- \hskip.3cm $(ab)^{\dag} =b^{\dag}a^{\dag}$,\par \indent
ii- \hskip.2cm $a^{\dag^*}(pq^{\dag}-q^{\dag}p)b=0$, and
$a^{\dag^*}(r^{\dag}s-sr^{\dag})b=0$,\par \indent iii-
$spq^{\dag}p=q^{\dag}p$, and $sr^{\dag}sp=sr^{\dag}$.\par
\end{theo3.4}
\begin{proof}
\indent The proof is similar to the one of Theorem 3.1. However,
instead of the third statement of Proposition 2.4, the fifth
statement of the mentioned proposition must be used. In addition, in
order to compute $q^{\dag}$ and $r^{\dag}$, Theorem 7 of [10] must be
considered.\par
\end{proof}

\newtheorem*{rem3.5}{Remark 3.5}
\begin{rem3.5}\rm Consider a $C^*$-algebra $A$, and two regular
elements of $A$, $a$ and $b$, such that $ab$ is also
regular. It is well-known that the reverse order law for the
product $ab$ is equivalent to the conditions
$$
a^{\dag}abb^*=bb^*a^{\dag}a,\hskip1cm  bb^{\dag}a^*a=a^*abb^{\dag},
$$
\noindent see for example Theorem 2 of [7], which was
proved for matrices but whose proof remains valid
in a $C^*$-algebra, Proposition 4.4
of [18], Remark 3.2 of [4],
Corollary 3.11 of [12], and also Theorem 5 of [13].
However, the above conditions are equivalent to
the equalities 
$$
pq=qp, \hskip1cm rs=sr.
$$
\indent In fact, the first condition is exactly
$$rs=sr.$$
\indent As regard the second condition, since $bb^{\dag}$
commutes with $a^*a$, according to Theorem 5 of [10],
$bb^{\dag}$ commutes with $(a^*a)^{\dag}$, which, according to
Theorem 7 of [10], proves that $pq=qp$.\par
\indent On the other hand, if $bb^{\dag}$
commutes with $a^{\dag}a^{\dag^*}$, according again to
Theorem 5 of [10], $bb^{\dag}$ commutes with the Moore-Penrose
inverse of $a^{\dag}a^{\dag^*}$. In particular, according
to Theorem 7 of [10], $bb^{\dag}$ commutes with $a^*a$.\par

\indent Furthermore, note that, according to Theorems 5 of [10], 
and to the fact that $p$, $q$, $r$ and $s$
are hermitian elements of $A$, the above conditions and
equalities are equivalent to
$$
q^{\dag}p=pq^{\dag},\hskip1cm r^{\dag}s=sr^{\dag}.
$$
\indent Consequently, the second condition of Theorems 3.1
and 3.2 (resp. Theorems 3.3 and 3.4) could have been replaced by
the commutativity of $p$ and $q$, and of $r$ and $s$
(resp. the commutativity of $q^{\dag}$ and $p$ and of $r^{\dag}$
and $s$), however, this has not been done for two reasons.
In first place, the conditions \it ii \rm in the aforesaid
Theorems are weakers, but above all, Theorems
3.1 - 3.4 have been presented in a way that they
provide a characterization of the reverse order law
for the Moore-Penrose inverse analogue to the one of Theorem
3.1 of [14] for generalized inverses.\par
\end{rem3.5}
\newpage
\bf{4. Moore-Penrose hermitian elements}\rm \vskip.5cm

\indent In first place, the main notion of this and the following section
is introduced.\par

\newtheorem*{def4.1}{Definition 4.1}
\begin{def4.1}\rm Consider a $C^*$-algebra $A$. A regular element
$a\in A$ will be called Moore-Penrose hermitian, if  $a^{\dag}=a$.
\par\end{def4.1}
\indent  Next follow the basic facts regarding the concept just introduced.
In the next section Moore-Penrose hermitian elements
will be fully characterize.\par

\newtheorem*{prop4.2}{Proposition 4.2}
\begin{prop4.2} Consider a $C^*$-algebra $A$ and an element $a\in A$.
Then, the following statements hold:\par i- \hskip.2cm Necessary and
sufficient for $a$ to be a Moore-Penrose hermitian element is
$a=a^3$ and $(a^2)^*=a^2$.\par ii- \hskip.1cm If $a$ is a
Moore-Penrose hermitian element, then $a^n$ also is, $n\in\Bbb
N$.\par iii- The element $a$ is Moore-Penrose hermitian if and only
if $a^*$is.\par iv- \hskip.1cm If $a$ is a Moore-Penrose hermitian
element, then $\sigma (a)\subseteq \{0, -1, 1\}$, where $\sigma (a)$
denotes the spectrum of $a$.\par
\end{prop4.2}
\begin{proof}
\indent Definition 4.1 and the equations defining the Moore-Penrose
inverse prove the first point, which in turn proves the second.\par
\indent The third point is clear, and the fouth is a consequence of
the fact that $p=a^2$ is a hermitian idempotent.\par
\end{proof}

\bf{5. Characterizations of Moore-Penrose hermitian elements}\rm
\vskip.5cm \indent This section begins with the characterization of
Moore-Penrose hermitian $C^*$-algebra elements. In first place, some
notation is given.\par
\par \indent Recall that if $A$ is a $C^*$-algebra and $a\in A$,
then $L_a\colon A\to A$ is the map defined by left multiplication by
$a$, that is
$$
L_a(x)=ax,\hskip2cm (x\in A).
$$
\indent In addition, the range and the null space of $L_a$ will be
denoted by $R(L_a)=aA$ and $N(L_a)=a^{-1}(0)$ respectively.\par
\newtheorem*{theo5.1}{Theorem 5.1}
\begin{theo5.1} Consider a $C^*$-algebra $A$. Then, the following
statements are equivalent:\par
\end{theo5.1} i- \hskip.2cm $a\in A$ is a Moore-Penrose hermitian element,\par
ii- \hskip.1cm $aA=a^*A$, $a^{-1}(0)= a^{*-1}(0)$, $A=aA\oplus
a^{-1}(0)$, and if $L=L_a\mid_{aA}\colon aA\to aA$ and
$\tilde{L}=L_{a^*}\mid_{aA}\colon aA\to aA$, then
$L^2=\tilde{L}^2=\tilde{I}$, where $\tilde{I}$ denotes the identity
map of $aA$.\par
\begin{proof}
\indent Suppose that $a$ is a Moore-Penrose hermitian element of
$A$, and consider the map $L_{a^2}\colon A\to A$. Since $a^2$ is an
idempotent, $L_{a^2}$ is a projection defined in $A$. Consequently,
$A=R(L_{a^2})\oplus N(L_{a^2})$. However, since $a$ is a
Moore-Penrose hermitian element of $A$, an easy calculation proves
that $R(L_{a^2})=aA$ and $N(L_{a^2})=a^{-1} (0)$. Moreover, since
$L_{a^2}$ is a projection, it is clear that $L^2=\tilde{I}$.\par
\indent In addition, according to the third statement of Proposition
4.2, $a^*$ is a Moore-Penrose hermitian element. Moreover, according
to the fifth statement of Proposition 2.4, $a=a^*a^*a$ and
$a^*=aaa^*$, which implies that $aA=a^*A$. Furthermore, since
according to the third statement of Proposition 2.4, $a=aa^*a^*$ and
$a^*=a^*aa$, it follows that $a^{*-1}(0)=a^{-1}(0)$. However, since
$a^2$ is hermitian, $\tilde{L}^2=L^2=\tilde{I}$.\par \indent
Conversely, if the second statement holds, a straightforward
calculation porves that $L_a=L_a^3$ and $L_{a^*}^2=L_a^2$, which
clearly implies that $a=a^3$ and $(a^2)^*=a^2$, that is $a$ is a
Moore-Penrose hermitian element.\par
\end{proof}

\indent Next follows the characterization of Moore-Penrose hermitian
ele-ments in the frame of Hilbert spaces. However, firstly several
notions and results need to be reviewed\par \indent As in the case
of a $C^*$-algebra element, a Hilbert space operator will be said
\it Moore-Penrose hermitian\rm, if it has a Moore-Penrose inverse
$T^{\dag}$ and
$$
T^{\dag}=T.
$$
 \indent Recall that if
$T\colon H\to H$ is a bounded linear operator defined on the Hilbert
space $H$, then the Moore-Penrose inverse of $T$ is the unique
linear and continuous map $T^{\dag}$ for which the following
equations hold:\par
$$
T=TT^{\dag}T, \hskip.3cm T^{\dag}=T^{\dag} TT^{\dag},\hskip.3cm
(TT^{\dag})^*=TT^{\dag},\hskip.3cm (T^{\dag}T)^*=T^{\dag}T.
$$
\indent Note that the operator $T$ admits a generalized inverse in
$A=L(H)$ if and only if $R(T)$ is closed, see Theorem 3.8.2 of [9].
However, when the operator $T$ admits a generalized inverse, it can
be chosen to be the Moore-Penrose inverse of $T$ in $A=L(H)$, see
Theorem 5 of [10]. Moreover, in this case it is unique, and it
coincides with the Moore-Penrose inverse of $T$ viewed as an
operator defined on $L(H)$, see [10], [11], [14] and [15]. In addition,
a bounded linear map which has a generalized inverse will be called
a \it regular operator\rm.\par

\newtheorem*{theo5.2}{Theorem 5.2}
\begin{theo5.2} Consider a Hilbert space $H$, and $T$ a regular
operator defined on $H$. Then the following statements are
equivalent:\par \indent i- \hskip.2cm $T$ is a Moore-Penrose
hermitian operator,\par \indent ii- \hskip.1cm there exist two
orthogonal Hilbert subspaces $H_1$ and $H_2$ such that $H=H_1\oplus
H_2$, $T\mid H_1=0$, $T(H_2)\subseteq H_2$, and if $T_2$ denotes the
restriction of $T$ to $H_2$, then $T_2^2=I_2$, where $I_2$ denotes
the identity map of $H_2$.\par
\end{theo5.2}
\begin{proof}Suppose that $T^{\dag}=T$ and consider the self-adjoint
projection $P=T^{\dag}T=TT^{\dag}=T^2$. In particular, the Hilbert
space can be pre-sented as the orthogonal direct sum
$H=N(T^2)\oplus R(T^2)$. However, as in the case of a $C^*$-algebra,
since $T$ is a Moore-Penrose hermitian operator, a straightforward
calculation proves that $N(T)=N(T^2)$ and $R(T)=R(T^2)$. Define
$H_1=N(T)$, $H_2=N(T)^{\perp}=R(T)$, where $N(T)^{\perp}$ denotes
the orthogonal subspace of $N(T)$. Then, it is clear that
$T(H_2)\subseteq H_2$, and $T_2^2=I_2$.\par \indent Conversely, it
the second statement of the theorem holds, it is clear that $T^3=T$
and $T^2$ is the orthogonal projection onto $H_2$, in particular
$T^2$ is an hermitian projection.\par
\end{proof}
\indent Next normal Moore-Penrose hermitian elements will be
considered. However, first of all some preparation is necessary.\par
\indent Given a $C^*$-algebra $A$, the \it conorm \rm of an element
$a\in A$ is defined by
$$ c (a)=inf\{\parallel ax\parallel\colon dist (x,
a^{-1}(0))=1, \hbox{ } x\in A\},
$$
see [11] and [14].\par \indent It is worth noticing that if $a$ is a
regular element, then
$$
c(a)=\frac{1}{\parallel a^{\dag}\parallel},
$$
see Proposition 1.3 of [14] and Theorem 2 of [11].\par

\indent Next consider a bounded linear operator $T\colon H\to H$,
where $H$ is a Hilbert space. Then, $T$ is said a \it partial
isometry\rm, if $T$ admits a Moore-Penrose inverse and
$T^{\dag}=T^*$, see [15] and Chapter 15 of [8]. In order to keep an
analogy with the Hilbert space case, an element $a$ of a
$C^*$-algebra $A$ will be called a \it partial isometry\rm, if $a$
is regular and $a^{\dag}=a^*$.\par \indent It is clear that if $a\in
A$ is a partial isometry, then $a^*a$ is a hermitian idempotent.
Conversely, consider $a\in A$ such that $a^*a$ is a hermitian
idempotent. Then, since each $C^*$-algebra has a faithful
representation in a Hilbert space, see for example Theorem
 7.10 of [6], according to problem 127, Chapter 15, of [8], a
straightforward calculation shows that $a$ is a partial isometry.
Furthermore, since $a$ is a partial isometry if and only if $a^*$
is, then necessary and sufficient for $a$ to be a partial isometry
is that $aa^*$ is a hermitian idempotent. See [1] where an
equivalent definition of the notion under consideration was
considered.\par \indent In the following proposition a
generalization of Corollary 3.2 of [15] will be proved. This result
will be central for the characterization of normal Moore-Penrose
hermitian elements.\par
\newtheorem*{prop5.3}{Proposition 5.3}
\begin{prop5.3}Consider a $C^*$-algebra $A$, and a non-zero regular element
$a\in A$. Then, necessary and sufficient for $a$ to be a partial
isometry is $c(a)=\parallel a\parallel =1$.\par
\end{prop5.3}
\begin{proof} Let $a\in A$ be a non-zero regular element, and consider,
according to Theorem 7.10 of [6], a Hilbert space $H$ and $\pi\colon
A\to L(H)$ a faithful representation of $A$. It is worth noticing
that in this case $\pi (a)\in L(H)$ is regular and $\pi
(a)^{\dag}=\pi (a^{\dag})$.\par \indent Suppose that $a$ is a
partial isometry. Then $\pi (a)\in L(H)$ also is a partial isometry.
Then, according to Corollary 3.2 of [15], $\parallel\pi(a)\parallel$
$=1$. In particular, $\parallel a\parallel =1$. Moreover, according
to Proposition 1.3 of [14],
$$c(a)=\frac{1}{\parallel a^*\parallel}=\frac{1}{\parallel
a\parallel}=1.$$ \indent Conversely, suppose that a regular element
$a\in A$ is such that $c(a)=$ $\parallel a\parallel =1$. Then,
$\parallel \pi (a)\parallel =1$, and according again to Proposition
1.3 of [14], $\parallel a^{\dag}\parallel =1$.
\par On the other hand, since $\pi (a^{\dag})=\pi (a)^{\dag}$,
according to Corollaries 2.3 and 3.2 of [15], $\pi (a)$ is a partial
isometry. However, since $\pi\colon A\to L(H)$ is a faithful
representation, $a$ is a partial isometry.\par\end{proof}

\newtheorem*{theo5.4}{Theorem 5.4}
\begin{theo5.4} Consider a $C^*$-algebra $A$. Then, an element
$a\in A$ is a normal Moore-Penrose hermitian element if and only if
$a$ is a hermitian partial isometry.\par\end{theo5.4}
\begin{proof}
\indent Suppose that $a\in A$ is a normal Moore-Penrose hermitian
ele- ment. Then, according to Theorem 2.9 of [6], to the fourth
statement of Proposition 4.2, and to Corollary 1.6 of [14], $c(a)=1$.
Moreover, since $a$ is a Moore-Penrose hermitian element, according to
Proposition 1.3 of [14] or to Theorem 2 of [11], $\parallel
a\parallel=1$. Consequently, according to Proposition 5.3, $a$ is
a partial isometry. However, $a=a^{\dag}=a^*$, that is $a$ is
hermitian.\par \indent The converse is clear.\par
\end{proof}

\indent Note that if $T\colon H\to H$ is a linear and continuous
Hilbert space map, then $T$ is a normal Moore-Penrose hermitian
operator if and only if the map $T_2$ in Theorem 5.2 is a hermitian
unitary operator.\par
\vskip.5cm
\noindent \bf{Acknowledgements.} \rm The author wishes to express his
indebtedness to the referees of this article, for their remarks and suggestions
led to an improvement of the present work.\par

\vskip1cm
\hskip.4cm Enrico Boasso\par
E-mail address: $\hbox{enrico}\_\hbox{odisseo@yahoo.it}$\par

\end{document}